\documentclass[10pt]{article}
 
\usepackage[latin1]{inputenc}
\usepackage{enumerate} 
\usepackage[pagebackref=true, colorlinks=true, linkcolor=blue, citecolor=blue, linkcolor=blue, anchorcolor=red, urlcolor=blue]{hyperref}
\usepackage{url, amsmath,fancyhdr,amssymb, amsthm, 
pstricks, epsf, lscape, ifpdf, graphicx, float}
\usepackage{algorithm}
\usepackage{color}
\usepackage{algpseudocode} 
\algnewcommand\algorithmicinput{\textbf{Input:}}
\algnewcommand\INPUT{\item[\algorithmicinput]}
\algnewcommand\algorithmicinitialization{\textbf{Initialization:}}
\algnewcommand\INITIALIZATION{\item[\algorithmicinitialization]}

\usepackage{amsmath,amssymb,tikz-cd}

\usepackage{tabularx}
\newcolumntype{C}{>{\centering\arraybackslash}X}
\usepackage{amsmath}
\makeatletter
\newcommand{\leqnomode}{\tagsleft@true}
\newcommand{\reqnomode}{\tagsleft@false}
\makeatother

\newcommand{\Newchange}[1]{#1}

\usepackage[margin=3cm]{geometry}

\mathchardef\mhyphen="2D

\numberwithin{equation}{section}

\newtheorem{Definition}{Definition}
\newtheorem{Convention}{Convention}

\newtheorem{Example}{Example}
\newtheorem{Proposition}{Proposition}
\newtheorem{Lemma}{Lemma}
\newtheorem{Theorem}{Theorem}
\newtheorem{Corollary}{Corollary}
\newtheorem{Remark}{Remark}
\newtheorem{Assumption}{Assumption}

\newcommand{\ba}{\begin{array}}  
\newcommand{\ena}{\end{array}}

\newcommand{\A}{\mathcal A}

\newcommand{\trace}{\operatorname{trace}}

\newcommand{\closure}{\operatorname{closure}}

\newcommand{\dir}{\operatorname{dir}}

\newcommand{\eps}{\epsilon} 
\newcommand{\bpx}{\begin{pmatrix}}
\newcommand{\epx}{\end{pmatrix}}
\newcommand{\bbx}{\begin{bmatrix}}
\newcommand{\ebx}{\end{bmatrix}}

\newcommand{\front}{\operatorname{front}}

\newcommand{\bdef}{\begin{Definition}} 
\newcommand{\commentout}[1]{}
\newcommand{\co}[1]{}

\newcommand{\coab}[1]{}

\newcommand{\ti}{\times}

\newcommand{\norm}[1]{\parallel \! #1 \! \parallel}
\newcommand{\sym}[1]{{\cal S}^{#1}}
\newcommand{\psd}[1]{{\cal S}_+^{#1}}
\newcommand{\rad}[1]{\mathbb{R}^{#1}}

\newcommand{\eref}[1]{(\ref{#1})}

\newcommand{\R}{ {\cal R} }

\newcommand{\M}{ {\cal M} }

\newcommand{\dist}{\operatorname{dist}}

\newcommand{\beq}{\begin{equation}}
\newcommand{\eeq}{\end{equation}}
\newcommand{\beqa}{\begin{eqnarray}}
\newcommand{\eeqa}{\end{eqnarray}}
\newcommand{\bac}{\begin{array}{ccccccccccc}}
\newcommand{\eac}{\end{array}}
\newcommand{\bprop}{\begin{Proposition}}
\newcommand{\eprop}{\end{Proposition}}

\newcommand{\beqast}{\begin{eqnarray*}}
\newcommand{\eeqast}{\end{eqnarray*}}
\newcommand{\benum}{\begin{enumerate}}
\newcommand{\eenum}{\end{enumerate}}
\newcommand{\bit}{\begin{itemize}}
\newcommand{\eit}{\end{itemize}}
\newcommand{\bth}{\begin{Theorem}}
\newcommand{\enth}{\end{Theorem}}
\newcommand{\ble}{\begin{Lemma}}
\newcommand{\ele}{\end{Lemma}}
\newcommand{\bex}{\begin{Example}}
\newcommand{\eex}{\end{Example}}
\newcommand{\bcor}{\begin{Corollary}}
\newcommand{\ecor}{\end{Corollary}}
\newcommand{\brem}{\begin{Remark}}
\newcommand{\erem}{\end{Remark}}
\newcommand{\bass}{\begin{Assumption}}
\newcommand{\eass}{\end{Assumption}}

\newcommand{\bsmx}{\begin{small} \begin{pmatrix}}
\newcommand{\esmx}{\end{pmatrix} \end{small}}

\title{\Large Characterizing bad semidefinite programs: normal forms and short proofs}  
\author{G\'{a}bor  Pataki\thanks{Department of Statistics and Operations Research, University of North Carolina at Chapel Hill} \hspace{1cm} 
}

\begin{document}

\maketitle 

\begin{abstract}

 Semidefinite programs (SDPs) --  some of the most useful and versatile   optimization problems of the last few decades -- are often pathological: the optimal values 
 of the primal and dual problems may differ and  may not be attained.  
 Such SDPs 
 are  both theoretically interesting and  often  impossible to solve;  yet,   the pathological SDPs in the literature look strikingly similar. 
 
 Based on our recent work \cite{Pataki:17} we characterize pathological semidefinite systems by certain {\em excluded matrices}, which are easy to spot in all published examples. 
 Our main tool is a normal (canonical) form of 
 semidefinite systems, which makes their pathological behavior easy to verify. 
 The normal form is constructed in a surprisingly simple fashion, using mostly 
 elementary row operations inherited from Gaussian elimination. The proofs are elementary and  can be followed by a reader at the advanced undergraduate level.  
 
 As a byproduct, we  
 show how to transform any linear map acting on symmetric matrices 
 into a normal form, which allows us to quickly  check whether the  image of the semidefinite cone under the map  is closed.  
 We can thus introduce readers to a fundamental issue in convex analysis: 
 the linear image of a closed convex set may not be closed, and often simple conditions are available to verify 
 the closedness,  or lack of it. 

\end{abstract}



{\em OR/MS subject classification:} Primary: convexity; secondary: programming-nonlinear-theory

\section{Introduction. Main results} \label{sect-intro}

Semidefinite programs (SDPs)  -- optimization problems with 
semidefinite matrix variables, 
 a linear objective, and linear constraints -- are some of the most  practical, widespread, and interesting
  optimization problems of the last three decades. They naturally generalize linear programs, and 
 appear in diverse areas such as 
combinatorial optimization, polynomial optimization, 
engineering, and economics.
They are covered in many surveys, see e.g. 
\cite{Todd:00}  and textbooks, see e.g. \cite{BorLewis:05, Barvinok:2002, Ren:01, BonnShap:00,  BoydVand:04,   BentalNem:01,   Laurent-Vallentin:2016, Tuncel:11}.

They are also 
a subject of intensive research: 
in the last  30 years several thousand papers have been published on SDPs. 

To ground our discussion, let us write an  SDP in the form 
\begin{equation} \label{sdp-p} \tag{\mbox{$\mathit{SDP \mhyphen P}$}}     
\begin{array}{rl} 
\sup &  \sum_{i=1}^m c_i x_i  \\
s.t.   & \sum_{i=1}^m x_i A_i \preceq B,  
\end{array}
\end{equation}
where  $A_1, \dots, A_m, \,$ and $B$ are $n \times n$ 
symmetric matrices, $c_1, \dots, c_m$ are scalars, 
and for symmetric matrices $S$ and $T,$  we write 
$S \preceq T$ to say that $T - S$ is positive semidefinite (psd).  

To solve (\ref{sdp-p}) we rely on a natural dual, namely 
\begin{equation} \label{sdp-d}  \tag{\mbox{$\mathit{SDP \mhyphen D}$}}   
\begin{array}{rl}
\inf & B \bullet Y \\
s.t. & A_i \bullet Y = c_i \, (i=1, \dots, m)  \\ 
      & Y \succeq 0, 
      \end{array}
      \end{equation}
 where the inner product of symmetric matrices $S$ and $T$ is 
 $S \bullet T := \trace(ST).  \,$ 
Since the weak duality inequality 
\begin{equation} \label{eqn-weakduality} 
\sum_{i=1}^m c_i x_i \leq B \bullet Y
\end{equation}
always holds between feasible solutions $x$ and $Y, \,$ 
if  a pair $(x^*, Y^*)$  
satisfies (\ref{eqn-weakduality}) with equality, then they 
are both optimal.
Indeed, SDP solvers seek to find such an $x^*$ and $Y^*.$

However, SDPs often behave pathologically: the optimal values of 
(\ref{sdp-p}) and (\ref{sdp-d}) may differ and may not be attained. 

  The duality theory of SDPs -- together  with their  pathological behaviors -- is covered 
  in several references on optimization theory and in 
  textbooks written for broader audiences. 
  For example, \cite{BorLewis:05} gives  an extensive, yet concise account of Fenchel duality;  \cite{Todd:00} and \cite{Ren:01} 
  provide  very succinct treatments;  \cite{Barvinok:2002}  treats SDP duality as special case of  duality theory in infinite dimensional spaces;  
   \cite{BonnShap:00} covers stability and sensitivity analysis; 
  \cite{BentalNem:01} and \cite{BoydVand:04} contain many engineering applications;
   \cite{Laurent-Vallentin:2016} and \cite{Tuncel:11} 
  are accessible to  an audience with combinatorics background;  and \cite{Blekhetal:12} explores connections to algebraic geometry. 
    
Why are the pathological behaviors interesting?
First, they do not appear in linear programs, which makes it apparent that 
SDPs are a much less innocent generalization of linear programs, than one may think at first.
Note that the pathologies can come in ``batches": in extreme cases  (\ref{sdp-p}) and
(\ref{sdp-d}) {\em both} can  have unattained, and different, optimal values!  
The variety of thought-provoking  pathological SDPs makes teaching  SDP duality 
(to students mostly used to clean and pathology-free linear programming)  
a truly rewarding experience. 

Second, these pathologies  also appear in other convex optimization problems, 
thus SDPs  make excellent ``model problems" to study. 

Last but not least: pathological SDPs are often  difficult or impossible to solve.

Our recent paper  \cite{Pataki:17} was 
motivated by the curious similarity of pathological SDPs   in the literature. 
To build intuition, we recall two examples; they or their variants appear  in a number of  papers
 and surveys. 

\begin{Example} \label{ex1}  In the SDP  
	\begin{equation} \label{ex1-problem}
	\begin{array}{rl}
		\sup &   2 x_1 \\
		s.t. & x_1 \bpx 0 & 1 \\ 1 & 0 \epx \preceq \bpx 1 & 0 \\ 0 & 0 \epx
	\end{array} 
	\end{equation}
\Newchange{any feasible solution must satisfy $\bigl( \begin{smallmatrix}1 &  - x_1 \\ - x_1  &   0 \end{smallmatrix}\bigr) \succeq 0, \,$ i.e., $ - x_1^2 \geq 0, \,$ 
	so the only feasible solution is $x_1 = 0.$ } 
	
	The dual, with a variable matrix $Y = (y_{ij})$, is equivalent to 
	\begin{equation} \label{ex1-problem-dual} 
		\begin{array}{rllll}
		\inf & y_{11} \\
		s.t. & \bpx y_{11} & 1 \\ 1 & y_{22} \epx \succeq 0, 
		\end{array}
		\end{equation}
	so it has an unattained $0$ infimum. 
\end{Example} 

\begin{figure}[htp]
		\centering
		\includegraphics[width = 5cm]{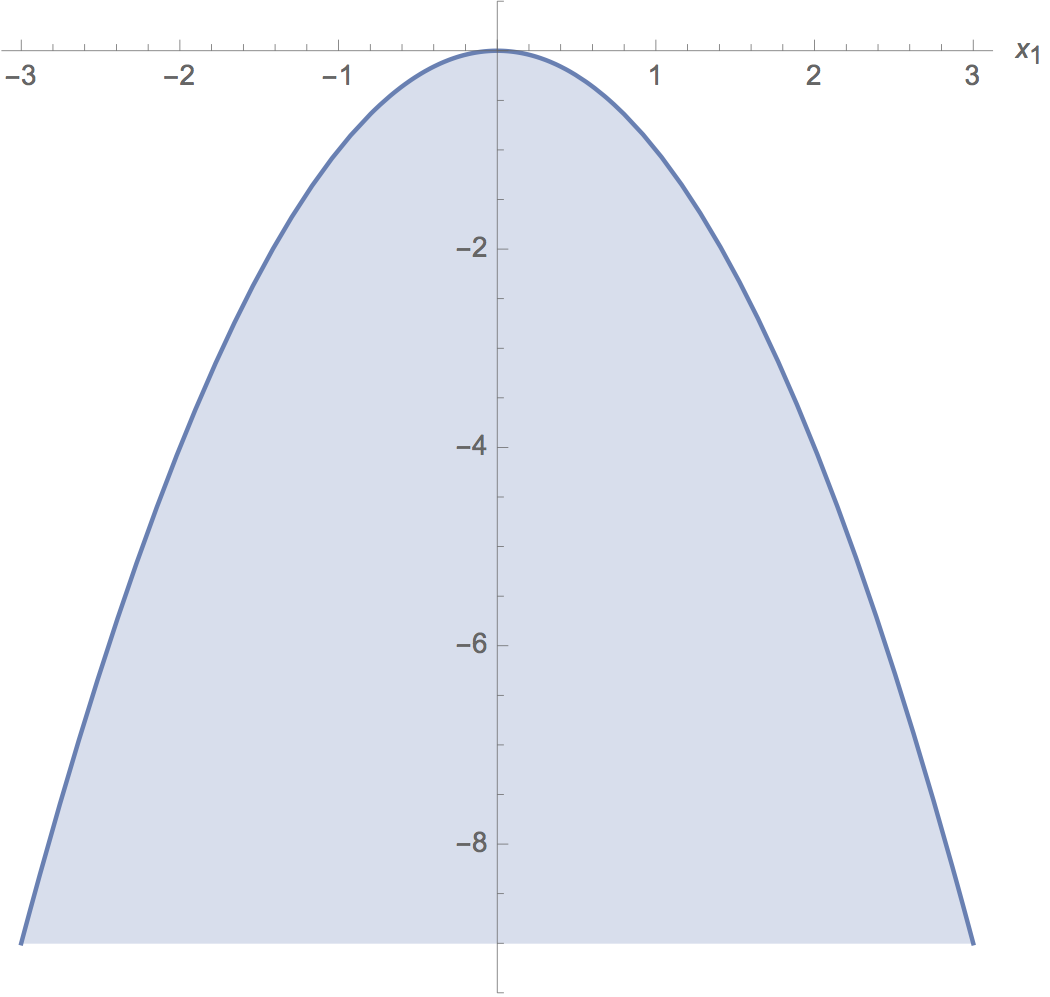} \hspace{1cm} \includegraphics[width = 5cm]{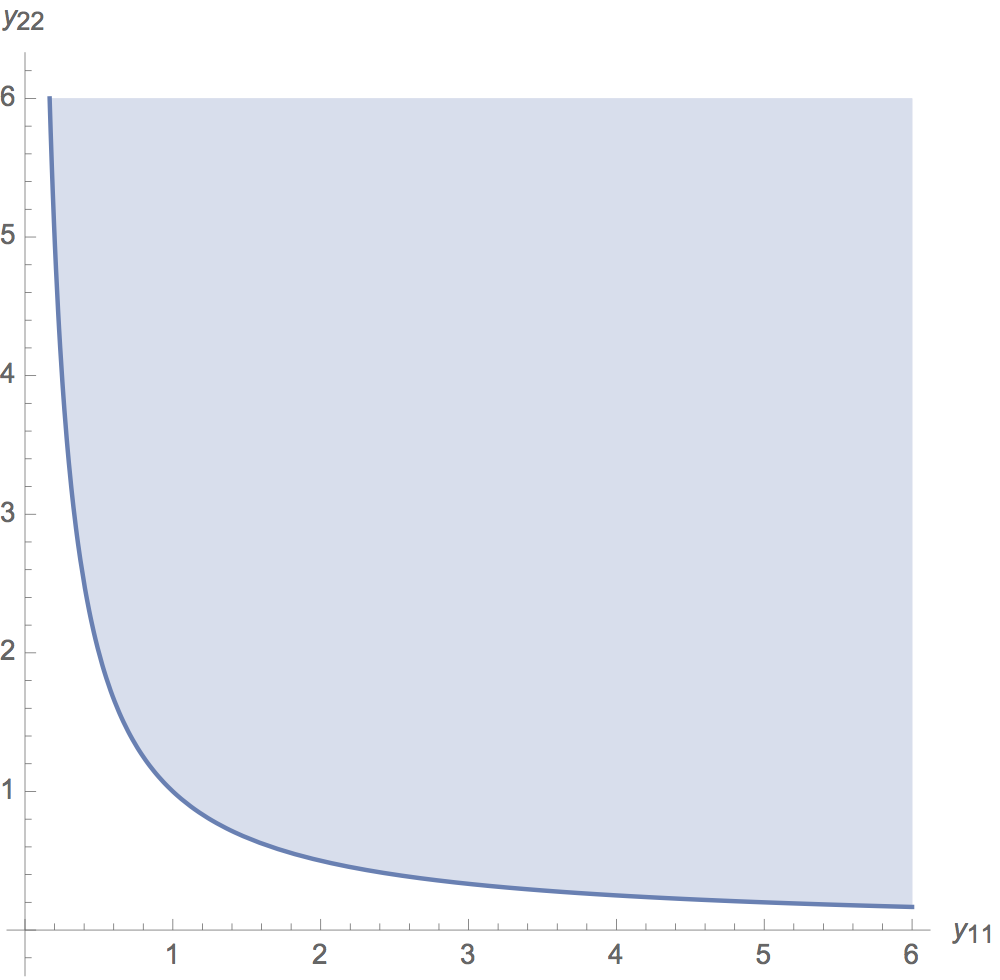}
		\caption{Parabola for the primal SDP, vs. hyperbola for the dual SDP in Example \ref{ex1}} 
		\label{figure-hyperbola} 
	\end{figure}

	Example \ref{ex1} has an  interesting connection to conic sections. 
	The primal SDP (\ref{ex1-problem}) seeks $x_1$ such that $- x_1^2 \geq 0, \,$ meaning 
	a point with nonnegative $y$-coordinate on a downward parabola. 
	This point is unique, so our parabola is ``degenerate."
	The dual 
	(\ref{ex1-problem-dual}) seeks the smallest nonnegative $y_{11}$ such that 
	$y_{11} y_{22} \geq 1, \,$ 
	i.e., the leftmost point on a hyperbola. This point, of course, does not exist: see Figure \ref{figure-hyperbola}. 
	
\begin{Example} \label{ex2}  We claim that the SDP 
	\begin{equation} \label{problem-ex2} 
	\begin{array}{rl}
	\sup  &  x_2   \\
	s.t. & x_1 \bpx 1 & 0 & 0 \\ 0 & 0 & 0 \\ 0 & 0 & 0 \epx + x_2 \bpx 0 & 0 & 1 \\ 0 & 1 & 0 \\ 1 & 0 & 0 \epx \preceq 
	\bpx 1 & 0 & 0 \\ 0 & 1 & 0 \\ 0 & 0 & 0 \epx 
	\end{array}
	\end{equation}
	has an optimal value that  differs from that of its dual. Indeed, in \eqref{problem-ex2} 
	\Newchange{we have $x_2 = 0$ in any feasible solution: this follows by a reasoning 
		 analogous to the one  we used in Example \ref{ex1}. Thus (\ref{problem-ex2}) has an attained $0$ supremum.} 
	
	 On the other hand, letting $Y = (y_{ij})$ be the dual 
		variable matrix,  the first dual constraint implies 
		$y_{11}=0.$ By $Y \succeq 0$ the first row and column of $Y$ is zero. By the second dual constraint $y_{22}=1$ so  the 
		 optimal value of the dual is $1,$ hence 
indeed there is a finite, positive duality gap. 
\end{Example}
Curiously, while their pathologies differ, 
Examples \ref{ex1} and \ref{ex2} still look similar. 
First, in both examples 
a matrix on the left hand side has a certain ``antidiagonal"  structure. Second, 
if we delete the second row and 
second column in all matrices in Example \ref{ex2}, and remove the first matrix, 
we get back Example \ref{ex1}! 
This raises the following questions:
Do all pathological semidefinite systems  ``look the same''?
Does the system of Example \ref{ex1} appear in all of them as a ``minor''? 

The paper \cite{Pataki:17} made these questions precise and gave a ``yes" answer to both.

To proceed, we state our main assumptions and recap needed  terminology from \cite{Pataki:17}. 
We assume throughout that (\ref{p-sd}) is feasible, and we say that the  semidefinite system 
\begin{equation} \label{p-sd} \tag{\mbox{$P_{\mathit{SD}}$}}
\sum_{i=1}^m x_i A_i \preceq B 
\end{equation}
is {\em badly behaved} if there is  $c \in \rad{m}$  for which 
the optimal value of (\ref{sdp-p})  is finite but the dual (\ref{sdp-d}) has no solution with the same value.  
We say that (\ref{p-sd}) is {\em well behaved}, if not badly behaved. 

A {\em slack matrix} or {\em slack}  in (\ref{p-sd})  is a psd  matrix of the form 
$ Z \, = \, B - \sum_{i=1}^m x_i A_i.  \,$ Of course, (\ref{p-sd}) has a 
maximum rank slack matrix, 
and our characterizations will rely on such a matrix. 

We also make the following assumption: 
\begin{Assumption} \label{ass-slack} 
	The maximum rank slack in (\ref{p-sd}) is 
	\begin{equation} \label{Zslack}
	Z \, = \, \bpx I_r & 0 \\ 0 & 0 \epx \, {\rm for \; some} \; 0 \leq r \leq n. 
	\end{equation}
\end{Assumption}
For the rest of the paper we fix this $r.$   

Assumption \ref{ass-slack} is easy to satisfy (at least in theory): 
	if $Z$ is a maximum rank slack in 
	(\ref{p-sd}),  and
	$Q$ is a matrix of suitably scaled eigenvectors of $Z, \,$ then replacing all $A_i$ by $Q^T A_i Q$ and $B$ by $Q^T B Q$ puts $Z$ into the required form.

A slightly strengthened version of the main result of 
\cite{Pataki:17} follows. 
\begin{Theorem}  \label{badsdp}
	The system \eref{p-sd} 
	is badly behaved if and only if the ``Bad condition" below holds:
	
	{\bf Bad condition:} There is a $V$ matrix, which is a linear combination of the 
	$A_i,$  and of the form 
	\begin{equation} \label{Vform}
	V \, = \, 
	\begin{pmatrix} V_{11} & V_{12} \\
	V_{12}^T & V_{22} 
	\end{pmatrix}, \, \, {\rm where} \, V_{11} \, {\rm is \,} r \times r, \, V_{22} \succeq 0, \, \R(V_{12}^T) \not \subseteq \R(V_{22}),
	\end{equation}
	where $\R()$ stands for rangespace.
	\qed \end{Theorem} 
The $Z$ and $V$ matrices are {\em certificates} of the bad behavior. 
They can be chosen as  
$$
Z \, = \, \bpx 1 & 0 \\ 0 & 0 \epx, \, V \, = \, \bpx 0 & 1 \\ 1 & 0 \epx \, \mbox{in Example \ref{ex1}},   \, \mbox{and} 
$$
$$
Z \, = \, \bpx 1 & 0 & 0 \\ 0 & 1 & 0 \\ 0 & 0 & 0 \epx, \, V \, = \, \bpx 0 & 0 & 1 \\ 0 & 1 & 0 \\ 1 & 0 & 0 \epx \, \mbox{in Example \ref{ex2}}. 
$$

Theorem \ref{badsdp}  is  appealing: it is simple, and the excluded matrices $Z$ and $V$ are easy to spot in essentially all badly behaved semidefinite systems in the literature. For instance, 
we invite the reader to spot 
 $Z$ and $V$ (after ensuring Assumption \ref{ass-slack})  
in the SDP
$$
\sup  \,\,  x_2 \,\,  s.t. \,\, \bpx x_2 - \alpha & 0 & 0 \\
                                              0                   & x_1 & x_2 \\
                                              0                   & x_2 & 0 \epx \, \preceq \, 0,
                                              $$
which is Example 5.79 in \cite{BonnShap:00}.  Here  $\alpha > 0$ is a parameter, and the gap between this SDP and its dual is $\alpha.$ 

More examples are  in  
\cite{Ramana:97, KlepSchw:12, VanBo:96, TunWolk:12, LuoSturmZhang:97, Tuncel:11};
e.g., in an example \cite[page 43]{Tuncel:11} {\em any} matrix  on the left hand side can serve as a $V$ certificate matrix! 
	Theorem \ref{badsdp} also easily certifies 
	the bad behavior of some SDPs coming from 
	polynomial optimization, e.g., of the SDPs in  \cite{waki2012strange}.

Theorem \ref{badsdp}  has  an interesting geometric interpretation. 
	Let $\dir(Z, \psd{n})$ be the set of {\em feasible directions} at $Z$ in $\psd{n}, \,$ i.e.,
	\begin{equation} \label{eqn-feasidir} 
	\dir(Z, \psd{n}) \, = \, \{ \, Y \, | \, Z + \epsilon Y \succeq 0 \, {\rm for \, some \,} \epsilon > 0 \, \}. 
	\end{equation}
	Then 
	$V$ is in the {\em closure} of $\dir(Z, \psd{n}), \,$ but 
	it is not a feasible direction (see \cite[Lemma 3]{Pataki:17}). 
	 That is, for small $\epsilon > 0$ the matrix 
	 $Z + \epsilon V$ is ``almost" psd, but not quite.  
	
\begin{figure}[H]
		\centering
		\includegraphics[width = 10cm, height=8cm]{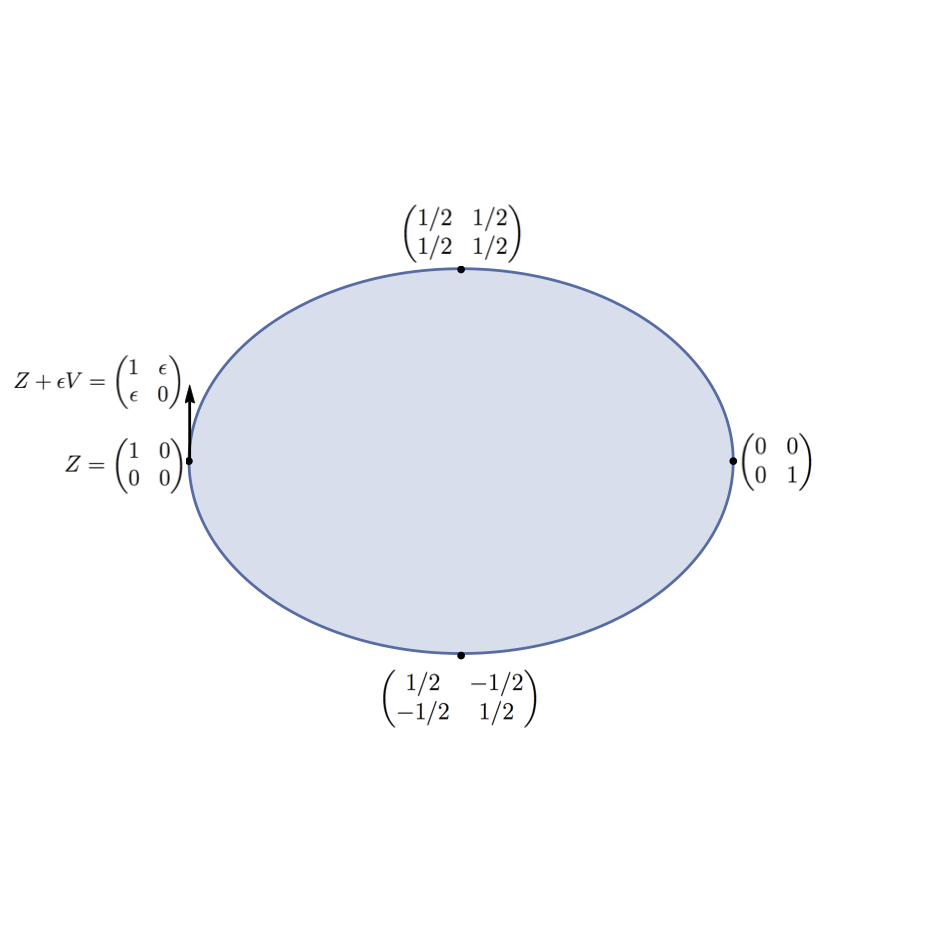} 
		\caption{The matrix $Z + \epsilon V$ is ``almost" psd, but not quite}
		\label{figure-ellipsoid}  
	\end{figure}

	We illustrate this point with the $Z$ and 
	$V$ of Example \ref{ex1}. The shaded region of Figure  
	\ref{figure-ellipsoid} 
	is the set of $2 \times 2$ psd matrices 
	with trace equal to $1.$ This set is an ellipse, so conic sections make a third appearance!
	The figure shows $Z$ and $Z + \epsilon V$ for a small $\epsilon >0.$

How do we characterize the good behavior of (\ref{p-sd})? We could, of course, say  that (\ref{p-sd}) is well behaved iff the $V$ matrix of Theorem \ref{badsdp} does {\em not} exist. However, 
there is a  much more convenient, and easier to check characterization, which we give below:

\begin{Theorem} \label{goodsdp}
	The system \eref{p-sd} is well behaved 
	if and only if both ``Good conditions" below hold.
	
	{\bf Good condition 1:} There is $U  \succ 0$ such that 
	\begin{equation}    \nonumber        
	A_i \bullet \bpx 0 & 0 \\ 0 & U \epx  \, = \, 0 \, {\rm for \, all \,} i. 
	\end{equation}

	{\bf Good condition 2:} 
	If $V$ is a  linear combination of the $A_i$ of the form 
	\begin{equation} \nonumber
	V \, = \, 
	\bpx V_{11}   & V_{12} \\
	V_{12}^T & 0  \epx, \,\, {\rm then} \, V_{12} = 0. 
	\end{equation}
	\qed
\end{Theorem} 
In Theorem \ref{goodsdp} and the rest of the paper, $U \succ 0$ means that $U$ is symmetric and positive definite, and we use the following convention: 
\begin{Convention}
  If  a  matrix is  partitioned  as in Theorems \ref{badsdp} or \ref{goodsdp}, then we understand that 
 the upper left block is $r \times r.$  
\end{Convention}

\begin{Example} \label{ex3} 
At first glance, the system 
	\begin{equation} \label{seinfeld} 
	x_1 \bpx 0 & 0 & 0 \\ 0 &  0 & 1 \\ 0 & 1 & 0 \epx \preceq \bpx 1 & 0 & 0 \\ 0 & 0 & 0 \\ 0 & 0 & 0 \epx
	\end{equation}
	looks very similar to the system in Example \ref{ex1}. 
	However, (\ref{seinfeld}) is well behaved, and Theorem \ref{goodsdp} verifies this by 
	choosing $U = I_2$ in 
	``Good condition 1"   (of course ``Good condition 2"  trivially holds). 
\end{Example}

In  \cite{Pataki:17} we proved  Theorems \ref{badsdp} and \ref{goodsdp}  
from a much more general result 
(Theorem 1 therein), which characterizes badly (and well) behaved conic linear systems.
In this paper we give short  proofs of Theorems \ref{badsdp} and \ref{goodsdp} using building blocks 
from \cite{Pataki:17}. Our proofs mostly use 
elementary linear algebra: we reformulate (\ref{p-sd}) into  normal  forms that make 
its bad or good behavior  trivial to recognize. The normal forms are  
inspired by  the row echelon form 
of a linear system of equations, and most of the operations that we use to construct them indeed come from 
Gaussian elimination.

As a byproduct,  we show how to construct 
normal forms of 
linear maps 
$$
\M: \, n \times n \, {\rm symmetric \, matrices} \, \rightarrow \rad{m}, 
$$
to easily verify whether the image of the cone of semidefinite matrices
under $\M$ is closed. We can thus introduce  students  
to a fundamental issue in convex analysis:
the linear image of a closed convex set is not always closed, and we can often verify 
its (non)closedness via simple conditions. For recent literature on closedness criteria see e.g.,   
\cite{BausBor:99,  auslender2006asymptotic, BertTseng:07, BorweinMoors:09, BorweinMoors:10, Pataki:07};  
for connections to duality theory, see  e.g. 
\cite[Theorem 7.2]{Barvinok:2002}, \cite[Theorem 2]{HenrionKorda:2014} , \cite[Lemma 2]{Pataki:17}.
\Newchange{For us the most relevant closedness criteria are in  \cite[Theorem 1]{Pataki:07}:  these  criteria 
led to the results of \cite{Pataki:17}.} 

We next describe how to {\em reformulate} (\ref{p-sd}). 
\begin{Definition} \label{definition-reform}  A semidefinite system is 
	an {\em elementary reformulation,}  or {\em reformulation}
	of (\ref{p-sd}) if it is obtained from (\ref{p-sd}) 
	by a  sequence of the following operations:
	\begin{enumerate}
		\item \label{rotate} Choose an invertible matrix of the form 
		$$T = \bpx I_r & 0 \\ 0 & M \epx,$$
		and replace $A_i$ by $T^T A_i T$ for all $i$ and $B$ by $T^T B T.$  
		\item \label{slack} Choose $\mu \in \rad{m}$ and replace $B$ by $B + \sum_{j=1}^m \mu_j A_j. \,$  
		\item \label{exch} Choose indices $i \neq j$ and exchange $A_i$ and $A_j.  \,$  
		\item \label{lambda} Choose $\lambda \in \rad{m}$ and an index $i$ such that $\lambda_i \neq 0, $ and replace 
		$A_i$ 
		by $\sum_{j=1}^m \lambda_j A_j.  \,$ 
	\end{enumerate}
\end{Definition}
(Of course, we can use just some of these operations and we can use them in any order). 

Where do these operations come from? As we mentioned above, mostly from Gaussian elimination: the last three 
can be viewed as elementary row operations done on  (\ref{sdp-d}) with some $c \in \rad{m}.$ 
For example,  operation (\ref{exch}) exchanges the constraints 
$$
A_i \bullet Y = c_i \,\, {\rm and} \,\, A_j \bullet Y = c_j.
$$
Reformulating (\ref{p-sd}) keeps the  maximum rank slack the same (cf. Assumption   \ref{ass-slack}).
Of course, (\ref{p-sd}) is badly behaved if and only if its reformulations are.

We organize the rest of the paper as follows. 
In the rest of this section we review preliminaries. 
In Section \ref{section-proofs} we prove Theorems \ref{badsdp} and \ref{goodsdp}
 and show how to construct the normal forms. 
We prove the chain of implications 
\begin{equation} \label{proofs-bad} 
\begin{array}{rcl} 
(\ref{p-sd}) \, \mbox{satisfies the ``Bad condition"} & \Longrightarrow & \mbox{it has a ``Bad reformulation"}    \\    
    & \Longrightarrow & \mbox{it is badly behaved}, 
\end{array}
\end{equation}
and the ``good" counterpart
\begin{equation} \label{proofs-good} 
\begin{array}{rcl} 
(\ref{p-sd}) \, \mbox{satisfies the ``Good conditions"} & {\Longrightarrow} & \mbox{it has a ``Good reformulation"}    \\  
& \Longrightarrow & \mbox{it is well behaved}. 
\end{array}
\end{equation}
In these proofs we only use elementary linear algebra. 

Of course, if (\ref{p-sd}) is badly behaved, then it is not well behaved. 
Thus the implication 
\begin{equation} \label{proofs-tie} 
\begin{array}{rcl} 
\mbox{Any of the ``Good conditions" fail}  & \Longrightarrow & \mbox{the ``Bad condition" holds},
\end{array}
\end{equation}
ties everything together and shows that in (\ref{proofs-bad}) and (\ref{proofs-good}) equivalence holds. 
Only the proof of (\ref{proofs-tie}) needs some elementary duality theory (all of which we recap in Subsection \ref{subsection-notation}), 
thus all proofs can be followed by a reader at the advanced undergraduate level.

In Section \ref{section-closed} we look at 
linear maps that act on 
symmetric matrices. As promised, we show how to bring them 
into a normal  form, 
 to easily check whether the image of the 
cone of semidefinite matrices under such a map is closed.   
We also point out connections to 
asymptotes of convex sets, and 
weak infeasibility in SDPs.
In Section \ref{section-conclusion} we close with a discussion.

\subsection{Notation and preliminaries} 
\label{subsection-notation}

As usual, we let  $\sym{n}$ be the set of $n \times n$ symmetric matrices,
and $\psd{n}$ the set of $n \times n$ symmetric positive semidefinite 
matrices. 

For completeness, we next prove the weak duality inequality (\ref{eqn-weakduality}). Let $x$ be feasible in 
(\ref{sdp-p}) 
and $Y$ be feasible in (\ref{sdp-d}). 
Then 
\begin{equation} \nonumber 
B \bullet Y - \sum_{i=1}^m c_i x_i \, = \, B \bullet Y - \sum_{i=1}^m (A_i \bullet Y) x_i \, = \, 
(B - \sum_{i=1}^m x_i A_i) \bullet Y \, \geq \, 0,
\end{equation}
where the last inequality follows, since the $\bullet$ product of two psd matrices is  nonnegative. Accordingly,  $x$ and $Y$ are both optimal iff the last inequality holds at equality.
 
 We next  discuss two well known regularity conditions, both of which   
 ensure that (\ref{p-sd}) is well behaved:  
\begin{itemize}
	\item The first is Slater's condition: this means that there is a positive definite slack in (\ref{p-sd}). 
	\item The second requires the $A_i$ and $B$ to be diagonal; in that case (\ref{p-sd}) is a polyhedron and (\ref{sdp-p}) is just a linear program. 
\end{itemize}
The sufficiency  of these 
conditions is immediate from Theorem \ref{badsdp}. If Slater's condition holds,  then 
$Z$ in Theorem \ref{badsdp} is just $I_n, \,$ so the $V$ certificate matrix cannot exist; if the 
$A_i$ and $B$ are diagonal, then so are their linear combinations, so again $V$  cannot exist. 

Thus Theorem \ref{badsdp} unifies these two (seemingly unrelated)   conditions, and 
we invite the reader to check that so does Theorem \ref{goodsdp}. 

We mention here that linear programs are sometimes also ``pathological," meaning  
 both primal  and dual may be infeasible. However, linear programs do not exhibit the pathologies that we study here.

\section{Proofs and examples} 
\label{section-proofs}

In this section we prove and illustrate the implications (\ref{proofs-bad}),  (\ref{proofs-good}),  and 
(\ref{proofs-tie}).
\subsection{The Bad} 
\label{subsection-bad} 

\subsubsection{From `` Bad condition" to  ``Bad reformulation"}  

We assume   the ``Bad condition" holds in \eref{p-sd} 
and show how to reformulate it as 
\begin{equation} \label{p-sd-bad} \tag{\mbox{$P_{SD, \mathit{bad}}$}} 
\sum_{i=1}^k x_i \bpx F_i & 0 \\ 0 & 0 \epx + \sum_{i=k+1}^m x_i \bpx F_i & G_i \\ G_i^T & H_i \epx \, \preceq \, \bpx I_r & 0 \\ 0 & 0 \epx = Z,
\end{equation}
where   
\begin{enumerate}
\item matrix $Z$ is the maximum rank slack, 
\item \label{indep}  matrices 
$$
\bpx G_i \\ H_i \epx \, (i=k+1, \dots, m) 
$$
are linearly independent, and 
\item $H_m \succeq 0.$
\end{enumerate}
Hereafter, we shall -- informally -- say that (\ref{p-sd-bad}) is a ``Bad reformulation" of (\ref{p-sd}).
We denote the constraint matrices on the left hand side by $A_i$ throughout the reformulation process.  

\Newchange{To begin,} we  replace $B$ by $Z$ in (\ref{p-sd}). 
We then choose $V = \sum_{i=1}^m  \lambda_i A_i$ to satisfy the ``Bad condition," 
and note that the block of $V$ comprising the last $n-r$ columns must be nonzero.
Next, we pick  an $i \,$ such that $\lambda_i \neq 0, \,$  and we use operation (\ref{lambda}) 
in Definition \ref{definition-reform} to  replace
$A_i$ by $V. \,$ We then 
switch $A_i$ and $A_m.$ 

Next we choose a maximal subset of the $A_i$ matrices whose 
blocks
comprising the last $n-r$ columns are linearly independent. We let $A_m$ be one 
of these matrices (we can do this since  $A_m$ is now the $V$ certificate matrix), 
and permute the $A_i$ so this special subset 
becomes $A_{k+1}, \dots, A_m$ for some $k \geq 0.$ 

\Newchange{Finally, we take linear combinations of the $A_i$}    to zero out the last $n-r$ columns of 
$A_1, \dots, A_k, \,$ and arrive at the required reformulation.
\qed 

Note that the systems in Examples \ref{ex1} and \ref{ex2} are already in the normal  form of (\ref{p-sd-bad}). 
The next example is  a counterpoint: it is  a more complicated badly behaved  system, which  at first 
is very far from 
being in the 	normal form.

\begin{Example}  \label{example-large-bad}  {\bf (Large bad example)}  
	The system 
	\begin{equation} \label{bad-orig} 
	\begin{split} 
x_1  \bpx 9  &   7 &    7 &    1 \\  7  &  12 &    8  &  -3 \\ 	7 &    8 &    2 &    4 \\ 	1  &  -3 &    4 &    0 \epx
+ x_2 \bpx 17  &  7 &   8 &   -1 \\ 7 &  8 &  7 & -3 \\ 8 &  7 &  4 &  2 \\ -1 & -3 & 2 & 0 \epx 
+ x_3 \bpx 1  &   2 &  2 &  1 \\ 2 &  6 & 3 & -1 \\ 2 &  3 &  0 &  2 \\ 1 &   -1 & 2 &  0 \epx \\ 
+ x_4 \bpx   9   & 6 &  7 &  1 \\ 	6 &  13 & 8 & -3 \\ 7 &  8 & 2 &  4 \\ 	1 & -3 & 4 &   0 \epx 
\preceq
\bpx 
45   & 26 &   29   &  2 \\
26   & 47 &   31 &  -12 \\
29 &   31 &   10 &   14 \\
2 &  -12 &   14 &    0 \epx 
\end{split}
\end{equation}
 is badly behaved, but this would be 
difficult to verify by any ad hoc method. 

\Newchange{Let us, however, verify its bad behavior using Theorem \ref{badsdp}.} 
System (\ref{bad-orig}) satisfies the ``Bad condition"  with 
$Z$ and $V$ certificate matrices 
\begin{equation} \label{Z-large-bad} 
Z \, = \,  \bpx 1 & 0 & 0 & 0 \\ 0 & 1 & 0 & 0 \\ 0 & 0 & 0 & 0 \\ 0 & 0 & 0 & 0 \epx, \, 
V = \bpx 7 &   2 &  3 & -1 \\ 2 &    1 &    2 &   -1 \\ 3 &    2 &  2 &  0 \\ -1 &  -1 &  0 & 0 \epx. 
\end{equation}
Indeed, $Z = B-A_1-A_2-2A_4, \, V = A_4 - 2 A_3 \,$ (where we write $A_i$ for the matrices on the left hand side, 
and $B$ for the right hand side), 
and we explain shortly why $Z$ is a maximum rank slack. 

Let us next reformulate  system (\ref{bad-orig}): after the  operations
\begin{equation} \label{operations} 
\begin{array}{rcl}
B & := & B - A_1 - A_2 - 2A_4, \\
A_4 & = & A_4 - 2A_3, \\
A_2 & = & A_2 - A_3 - 2A_4, \\
A_1 & = & A_1 - 2A_3 - A_4
\end{array}
\end{equation}
it becomes 
\begin{equation} \label{bad-reform} 
\begin{split}
x_1 \bpx  0 & 1 & 0 & 0 \\ 1 & -1 & 0 & 0 \\ 0 & 0 & 0 & 0 \\ 0 & 0 & 0 & 0 \epx 
+ x_2 \bpx 2 & 1 & 0 & 0 \\ 1 & 0 & 0 & 0 \\ 0 & 0 & 0 & 0 \\  0 & 0 & 0 & 0 \epx 
+ x_3 \bpx 1 & 2 & 2 & 1 \\ 2 & 6 & 3 & -1 \\  2 & 3 & 0 & 2 \\ 1 & -1 & 2 & 0 \epx \\
\hspace{2cm} + x_4 \bpx 7 & 2 & 3 & -1 \\ 2 & 1 & 2 & -1 \\ 3 & 2 & 2 & 0 \\ -1 & -1 & 0 & 0 \epx \preceq
   \bpx 1 & 0 & 0 & 0 \\ 0 & 1 & 0 & 0 \\ 0 & 0 & 0 & 0 \\ 0 & 0 & 0 & 0 \epx, 
   \end{split}
 \end{equation}
which is in the normal form of (\ref{p-sd-bad}). 
	Besides looking simpler than (\ref{bad-orig}), the bad behavior of (\ref{bad-reform}) is much easier  to verify, as we shall see soon.

How do we convince a ``user" that  $Z$ in  equation 
(\ref{Z-large-bad}) is indeed a maximum rank slack in  system (\ref{bad-orig}) ? 
Matrices 
\begin{equation} \label{y1y2}
Y_1 \, = \,  \bpx 0  &   0 &    0  &   0 \\
0  &   0 &    0  &   0 \\
0  &   0 &    0  &   0 \\
0  &   0 &    0  &   1 
\epx \, {\rm and} \, Y_2 \, = \,  \bpx 0  &   0 &    0  &   1 \\
0  &   0 &    0  &   1 \\
0  &   0 &    2  &   0 \\
1  &   1 &    0  &   0 
\epx
\end{equation}
have zero $\bullet$ product with all constraint matrices, and hence also with any slack.
\Newchange{Thus, if $S$ is any slack, then $S \bullet Y_1 = 0, \,$ so   
the $(4,4)$ element of $S$ is zero, hence the  entire 4th row and column of $S$ is zero (since $S \succeq 0$). 
 Similarly, 
$S \bullet Y_2 = 0$ shows  the 3rd row and column of $S$ is zero, thus the rank of $S$ is at most two. Hence
$Z$ indeed has maximum rank. 
} 
\end{Example}

In fact, Lemma 5 in \cite{Pataki:17} proves that (\ref{p-sd}) can always be reformulated, so that a similar sequence of matrices certifies that $Z$ has maximal rank.  To do so, we need to use 
operation (\ref{rotate})  in Definition \ref{definition-reform}.

\subsubsection{If (\ref{p-sd}) has a ``Bad reformulation,"   then it is badly behaved} 
\label{subsubsect-badref-badbehave} 

For this implication  we show that a system in the normal form of (\ref{p-sd-bad}) is badly behaved; and for that, we 
devise a simple  objective function which has a finite optimal value over 
(\ref{p-sd-bad}), while the dual SDP has no solution with the same value.  

To start, let  
$x$ be feasible in  (\ref{p-sd-bad})   with a corresponding slack $S.$ 
Observe that the last $n-r$ rows and columns 
of $S$ must be zero, otherwise $\frac{1}{2}(S+Z)$ would be
a slack with larger rank than $Z.$ 
Hence, by condition (\ref{indep}) (after the statement of (\ref{p-sd-bad})),  
we deduce
$x_{k+1} \, = \, \dots \, = \, x_m \, = \, 0,$ so the optimal value of the SDP
\begin{equation} \label{mofo}
\sup \, \{ \, - x_m \, | \, x \, \text{is feasible in } \mbox{\eref{p-sd-bad}} \, \}
\end{equation}
is  $0.$ 
We prove that its dual cannot have a feasible solution with value $0, \,$ 
so suppose that 
$$
Y \, = \, \bpx Y_{11} & Y_{12} \\ Y_{12}^T & Y_{22} \epx \succeq 0
$$
is such a solution. By $Y \bullet Z = 0$ we get $Y_{11} = 0, \,$ and since $Y \succeq 0$ 
we deduce $Y_{12} = 0.$ 
Thus 
$$
\bpx F_m & G_m \\ G_m^T & H_m \epx \bullet Y = H_m \bullet Y_{22} \geq 0,
$$
so $Y$ cannot be feasible in the dual of (\ref{mofo}), a contradiction. 
\qed

\begin{Example} (Example \ref{example-large-bad} continued)  
	Revisiting this example, the bad behavior 
	of (\ref{bad-orig}) is  nontrivial to prove, whereas that of (\ref{bad-reform}) 
	is easy: the objective function $\sup - x_4$ gives a $0$ optimal value over it, 
	while there is no dual solution with the same value. 
\end{Example}

\subsection{The Good} \label{subsection-good} 
\subsubsection{From ``Good conditions" to  ``Good reformulation"}

Let us assume that both "Good conditions" hold. We show how to reformulate (\ref{p-sd}) as 
\begin{equation} \label{p-sd-good} \tag{\mbox{$P_{SD, \mathit{good}}$}} 
\sum_{i=1}^k x_i \bpx F_i & 0 \\ 0 & 0 \epx + \sum_{i=k+1}^m x_i \bpx F_i & G_i \\ G_i^T & H_i \epx \preceq \bpx I_r & 0 \\ 0 & 0 \epx = Z,
\end{equation}
with the following attributes: 
\begin{enumerate}
\item \label{cond-good-1}  matrix $Z$ is the maximum rank slack.
\item \label{cond-good-2}  matrices 
$
H_i \, (i=k+1, \dots, m) 
$
are linearly independent.
\item \label{cond-good-3} $H_{k+1} \bullet U = \dots = H_m \bullet U = 0$ for some $U \succ 0.$ 
\end{enumerate}

We shall -- again informally -- say that (\ref{p-sd-good}) is a ``Good reformulation" of (\ref{p-sd}).
We construct the system (\ref{p-sd-good}) quite similarly to 
how we constructed (\ref{p-sd-bad}), and, as usual, we denote the matrices on the left hand side by $A_i$ throughout the process. 

We first replace $B$ by $Z$ in (\ref{p-sd}). 
We then  
choose a maximal subset of the $A_i$ whose 
lower principal $(n-r) \times  (n-r)$ blocks are linearly independent, and  
permute the $A_i,$ if needed, to make this subset 
$A_{k+1}, \dots, A_m$ for some $k \geq 0.$

Finally we take linear combinations to zero 
out the lower principal  $(n-r) \ti (n-r)$ block of $A_1, \dots, A_k.$ 
 By ``Good condition 2" 
the upper right $r \ti (n-r)$  block of 
$A_1, \dots, A_k$ 
(and the symmetric counterpart)
 also become zero. Thus items (\ref{cond-good-1}) and (\ref{cond-good-2}) hold. 
 
 As to item (\ref{cond-good-3}), 
 suppose  $U \succ 0$ satisfies ``Good condition 1." Then 
 $U$ has zero $\bullet$ product with the lower principal
  $(n-r)  \times (n-r)$ 
 blocks of the $A_i,$ hence $H_i \bullet U = 0$ for $i=k+1, \dots, m.$ 
Hence item  (\ref{cond-good-3}) holds, 
 and the proof is complete. 
\qed

\begin{Example} \label{example-large-good}
	  {\bf (Large good example)}  
	The system 
	\begin{equation} \label{good-orig} 
	\begin{split}
	x_1  \bpx 9  &  7 &    7 &  1 \\ 	7 &   12 &  8 &  -3 \\ 	7 &  8 &    2 &    4 \\ 1 &   -3 &  4 & -2 \epx
	+ x_2 \bpx 17  &   7 &   8 &   -1 \\ 	7 &    8 &  7 &   -3 \\ 8 &    7 &    4 &    2 \\ 	-1 &   -3 &    2 &   -4	 \epx 
	+ x_3 \bpx 1  &   2 &    2 &    1 \\ 	2 &   6 &  3 & -1 \\ 2 &    3 &    0 &    2 \\ 	1 &   -1 &   2 &  0 \epx \\
	+ x_4 \bpx  9  &   6  &   7 &   1 \\ 6 &   13 &    8 &   -3 \\ 	7 &    8 &  2 &   4 \\ 	1 &   -3 &    4 &   -2   \epx 
	\preceq
	\bpx 
	45  &  26 &   29 &    2 \\
	26  &  47  &  31  & -12 \\
	29 &   31  &  10  &  14 \\
	2 &  -12 &   14 &  -10
	\epx 
	\end{split}
	\end{equation}
	is well behaved, but it would be difficult to 
	improvise a method to verify this. 
	
	Instead, let us  check that the ``Good conditions" hold: to do so, 
	we write $A_i$ for the matrices on the left, and $B$ for the right hand side. 
	
	First, we can see that  "Good condition 1" holds with  $U = I_2, \,$ since  
	$$
	Y := \bpx 0 & 0 & 0 & 0 \\  0 & 0 & 0 & 0 \\  0 & 0 & 1 & 0 \\  0 & 0 & 0 & 1 \epx
	$$
	has zero $\bullet$  
	product with all $A_i$ (and also with $B$).  Luckily, $Y$ also certifies that 
   $Z$ in equation 
	(\ref{Z-large-bad}) is  a  maximum rank slack in (\ref{good-orig}): as $Y$ has 
	zero $\bullet$ product with any slack, the rank of any slack is at most two. Of course, $Z$ is a 
	rank two slack itself, since 	$Z = B - A_1 - A_2 - 2A_4. \,$ 	
	
	\co{
	Luckily, $Y$ also certfies that 
	$Z$ in  equation 
	(\ref{Z-large-bad}) is  a  maximum rank slack in the system (\ref{good-orig}): since 
	$B \bullet Y$ is also zero, $Y$ has 
	zero $\bullet$ product with any slack, thus the rank of any slack is at most two. Of course, $Z$ is a 
	rank two slack itself, since 	$Z = B - A_1 - A_2 - 2A_4. \,$ 	
}

	Next, let us verify ``Good condition 2." 
		Suppose the lower right  $2 \times 2$ block of $V := \sum_{i=1}^4 \lambda_i A_i$ is zero. 
		Then by a direct calculation  
		$\lambda \in \rad{4}$ is a linear combination of vectors 
		$$
		(-2, 1, 3, 0)^T \, {\rm and} \, (1,0,0,-1)^T, 
		$$
		so the upper right $2   \times 2 $ block of $V$ (and its symmetric counterpart) is also zero, so ``Good condition 2" holds. 
	
	Now, the same operations that are listed  in equation (\ref{operations})  turn system (\ref{good-orig}) 
	into 
		\begin{equation} \label{good-reform}  
	\begin{split}
	x_1 \bpx  0 &    1 &    0 &    0 \\
	1 &   -1 &    0  &   0 \\
	0 &    0 &    0 &    0 \\
	0 &    0 &    0 &    0
	\epx 
	+ x_2 \bpx 2  &   1 &    0 &   0 \\
	1  &   0 &    0 &    0 \\
	0 &    0 &    0 &    0 \\
	0  &   0 &    0  &   0
	\epx 
	+ x_3 \bpx 
	1   &  2 &    2 &    1 \\
	2 &    6 &    3 &   -1 \\
	2 &    3  &   0 &    2 \\
	1  &  -1  &   2 &    0
	\epx \\
	+ x_4 \bpx 
	7  &   2 &    3 &   -1 \\
	2  &   1 &    2  &  -1 \\
	3  &   2 &    2  &   0 \\
	-1 &   -1 &    0 &   -2
	\epx \preceq
	\bpx 1 & 0 & 0 & 0 \\ 0 & 1 & 0 & 0 \\ 0 & 0 & 0 & 0 \\ 0 & 0 & 0 & 0 \epx,
	\end{split}
	\end{equation}
	which is in the normal form of (\ref{p-sd-good}). As we shall see soon, the good behavior of 
	(\ref{good-reform}) is much easier  to verify. 
	\end{Example}

\subsubsection{If (\ref{p-sd})  has a ``Good reformulation,"  then it is well behaved} 
\label{subsub-goodref-to-wellbehaved} 

For this implication we show that the system  (\ref{p-sd-good}) is well behaved; 
and for that, we let  $c$ be such that  
\begin{equation} \label{origgood} 
v := \sup \, \biggl\{  \, \sum_{i=1}^m c_i x_i \, | \,  x \; {\rm \; is \; feasible \; in \; } \mbox{(\ref{p-sd-good})}   \biggr\}  
\end{equation}
is finite. 
An argument like the one in Subsubsection \ref{subsubsect-badref-badbehave} proves that 
$
x_{k+1} = \dots = x_m = 0 
$
holds for any  $x$ feasible in \eref{origgood}, so 
\begin{equation} \label{reduced}
v \, = \, \sup \, \{ \, \sum_{i=1}^k c_i x_i \, | \,  \sum_{i=1}^k x_i F_i \preceq I_r \, \}.
\end{equation}
Since \eref{reduced} satisfies Slater's condition,   
there is $Y_{11}$ feasible in its dual with $Y_{11} \bullet I_r = v.$ 
 
We next choose a $Y_{22}$ symmetric matrix (which  may not be not positive semidefinite), 
such that 
$$  
Y := \bpx Y_{11} & 0 \\ 0 & Y_{22} \epx
$$
satisfies the equality constraints of the dual of \eref{origgood} 
(this can be done, by condition (\ref{cond-good-2})). 
We then replace $Y_{22}$ by $Y_{22} + \lambda U$ for some $\lambda > 0$ to make it psd: we can do this by a simple linesearch.
After this,   $Y$ is feasible in the dual of \eref{origgood} 
(by condition (\ref{cond-good-3})),  
and clearly $Y \bullet Z = v$ holds. 
The proof is now complete. 
\qed

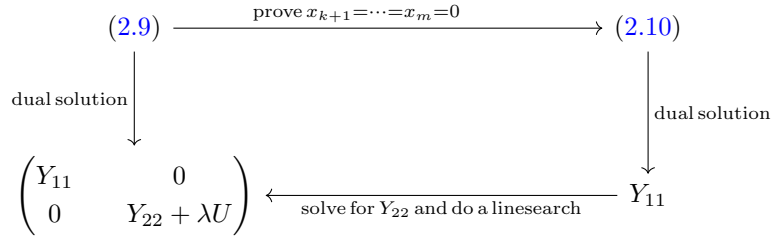
\begin{figure}[H]
		\begin{equation*} \label{eqn-commutative}
		\begin{tikzcd}[row sep=huge, column sep=width("bbbbbbbbbbbbbbbbbbbbbbb")]  
		(\ref{origgood})   \arrow[r,  "\qquad \qquad \qquad \qquad \quad  \mathrm{prove \, } x_{k+1} = \dots = x_m = 0 \hspace{4cm}"] 
		\arrow[d,   "\mathrm{dual \, solution}"  swap]  &     \arrow[d,  "\mathrm{dual \, solution}"] (\ref{reduced}) \\
		\bpx Y_{11} & 0 \\ 0 & Y_{22} + \lambda U \epx       & Y_{11}   \arrow[l, "\mathrm{solve \, for \,} Y_{22} \,  \mathrm{and \, do \,a \, linesearch}"] \\
		\end{tikzcd} 
		\end{equation*} 
		\vspace{-1.5cm} 
		\caption{How to construct an optimal  dual solution of (\ref{origgood})} 
		\label{figure-commutative} 
	\end{figure}

The above proof 
is illustrated in Figure 
\ref{figure-commutative} by a commutative diagram.
The horizontal arrows represent ``elementary" constructions, i.e., 
we find the object at the head of the arrow from the object at the tail of the arrow by  a basic argument or 
computation.    

\begin{Example} (Example \ref{example-large-good} (Large good example) continued.)   
We now  illustrate how to verify the good behavior of system 
	 (\ref{good-reform}): we pick an objective function with a finite optimal value over this system, and show how to construct an optimal dual solution.

We thus  consider the SDP 
\begin{equation} \label{problem-prove-good} 
\begin{array}{rl}
\sup & 2 x_2 + 5 x_3 + 7 x_4 \\
s.t.   & (x_1, x_2, x_3, x_4) \, {\rm is \, feasible \, in \, (\ref{good-reform})}, 
	\end{array}
	\end{equation}
in which $x_3 = x_4 = 0$ holds whenever $x$ is feasible, since in (\ref{good-reform}) the right hand side is the maximum 
	rank slack, and the lower right  $ 2 \times 2$ blocks of $A_3$ and $A_4$ are linearly independent. 

So the optimal value of (\ref{problem-prove-good}) is the same as that of 
\begin{equation} \label{problem-prove-good-reduced} 
\begin{array}{rl}
\sup & 2 x_2 \\
s.t.   & x_1 \bpx 0 & 1 \\ 1 & -1 \epx + x_2 \bpx 2 & 1 \\ 1 &  0 \epx \preceq \bpx 1 & 0 \\ 0 & 1 \epx.
\end{array}
\end{equation} 
\co{Next, let 
$$
(x_1, x_2) := \frac{1}{2} (-1,1), \, Y_{11} := \bpx 1 & 0 \\ 0 & 0 \epx, \, Y_{22} := \bpx 0 & 1 \\ 1 & 0 \epx, \, Y := \bpx Y_{11} & 0 \\ 0 & Y_{22} \epx.
$$
Here $(x_1, x_2)$ is  an optimal solution of (\ref{problem-prove-good-reduced}) and $Y_{11}$ is an optimal solution of its  dual: we can verify this directly. Further,  
$Y_{22}$ is chosen so that $Y$ satisfies the equality constraints of the dual of (\ref{problem-prove-good}). 
Next, let 
$$
 Y_{11} := \bpx 1 & 0 \\ 0 & 0 \epx, \, Y_{22} := \bpx 0 & 1 \\ 1 & 0 \epx, \, Y := \bpx Y_{11} & 0 \\ 0 & Y_{22} \epx.
$$
Here $(x_1, x_2)$ is  an optimal solution of (\ref{problem-prove-good-reduced}) and $Y_{11}$ is an optimal solution of its  dual: we can verify this directly. Further,  
$Y_{22}$ is chosen so that $Y$ satisfies the equality constraints of the dual of (\ref{problem-prove-good}). 
}
Next, let 
$$
Y_{11} := \bpx 1 & 0 \\ 0 & 0 \epx, \, Y_{22} := \bpx 0 & 1 \\ 1 & 0 \epx, \, Y := \bpx Y_{11} & 0 \\ 0 & Y_{22} \epx.
$$
Here  $Y_{11}$ is an optimal solution of the dual of (\ref{problem-prove-good-reduced}): this follows since it has the same value as the 
primal optimal solution $(x_1, x_2) = (-\frac{1}{2}, \frac{1}{2}).$ Further,  
$Y_{22}$ is chosen so that $Y$ satisfies the equality constraints of the dual of (\ref{problem-prove-good}).

	Of course, $Y_{22}$ is not psd, hence neither is $Y.$ As a remedy, we replace $Y_{22}$ by $Y_{22} + \lambda I_2$ 
	for some $\lambda \geq 1.$ This operation makes $Y$ feasible, 
	because  $U := I_2$ verifies item (\ref{cond-good-3}) (after the statement of (\ref{p-sd-good})).  
	Now  $Y$ is optimal in the dual of (\ref{problem-prove-good}) 
	and the process is complete. 
	
\end{Example} 

We remark that the procedure of constructing $Y$ from $Y_{11}$ was recently generalized in 
\cite{permenter2014partial} to the case when 
(\ref{p-sd}) satisfies only ``Good condition 2."

\subsection{Tying everything together} 
\label{subsection-tye-together} 

Now 
we tie everything together: 
we show that if any of the ``Good conditions" fail, then the ``Bad condition" holds.

Clearly, if  ``Good condition 2" fails, then the ``Bad condition" holds, 
so assume that ``Good condition 1" fails. 

First, we shall produce a matrix $V$ which is a linear combination of the $A_i$ such that 
\begin{equation} \label{eqn-V} 
V \, = \, \bpx V_{11} & V_{12} \\ V_{12}^T & V_{22} \epx \, {\rm with} \, V_{22} \succeq 0, \, V_{22} \neq 0. 
\end{equation} 
To achieve that goal, we let $B_i$ be the lower right order $n-r$ principal 
block of $A_i$ for $i=1, \dots, m$ and for some $\ell \geq 1$ 
choose matrices 
$C_1, \dots, C_\ell$ such that the set of their linear combinations is 
$$
 \{ \, U \in \sym{n-r}: B_1 \bullet U = \dots = \dots = B_m \bullet U = 0 \, \}.
$$
Consider next the primal-dual pair of SDPs 
\vspace{-.5cm} 
\begin{center}
	\begin{minipage}{0.5\linewidth}
		\leqnomode
		\begin{equation}\label{redp}
		\begin{split}
		\sup & \,\, t  \\
		s.t.   & \,\,  t I + \sum_{i=1}^\ell   x_i C_i \preceq 0 \\
		\end{split} 
		\end{equation}
	\end{minipage} \hspace{-2cm} 
	\begin{minipage}{0.5\linewidth}
		\begin{equation}\label{redd}
		\begin{split}
		\inf  & \,\, 0   \\
		s.t. & \,\, I \bullet W = 1 \\
		      & \,\, C_i \bullet W = 0  \, (i=1, \dots, \ell) \\
		& \,\, W \succeq 0. 
		\end{split} 
		\end{equation}
	\end{minipage}
\end{center}
Since ``Good condition 1" fails, the primal (\ref{redp}) has optimal value zero. 
The primal (\ref{redp})  also satisfies Slater's condition (with $x=0$ and $t=-1$) so the dual (\ref{redd}) has a feasible solution $W.$ This $W$ is of course nonzero, and a linear combination of 
the $B_i,$ say 
$$
W = \sum_{i=1}^m  \lambda_i B_i \,\, {\rm for \, some \,} \lambda \in \rad{m}. 
$$
Thus, $V := \sum_{i=1}^m   \lambda_i A_i$ passes requirement (\ref{eqn-V}). 

We are done if we show $\R(V_{12}^T) \not \subseteq \R(V_{22}), \,$ so 
 assume otherwise,  i.e.,  assume 
$V_{12}^T = V_{22}D$ for some \mbox{$D \in \rad{(n-r) \times r}.$}  
Define 
$$
M = \bpx I & 0 \\ -D & I \epx,
$$
and replace $A_i$ by $M^T A_i M$ for all $i$ and $B$ by $M^T B M.$ 
After this, 
the maximum rank slack $Z$ in (\ref{p-sd}) remains the same (see equation (\ref{Zslack})) and 
$V$ is transformed into 
$$
M^TVM = \bpx V_{11} - D^T V_{12}^T    & 0  \\ 0 &  V_{22} \epx.
$$
Since $V_{22} \neq 0, \,$ we deduce  $Z+ \epsilon V$ has larger rank than $Z$ for a small 
$\epsilon >0, \,$ which is a contradiction. The proof is complete.
\qed

We thus proved the following corollary:

\begin{Corollary} \label{corollary-badlywell} 
	The system (\ref{p-sd}) is badly behaved if and only if it has a bad reformulation of the  form (\ref{p-sd-bad}). 
	
	It is well behaved if and only if  it has a good reformulation  of the  form (\ref{p-sd-good}). 
	
\end{Corollary}

 
 \brem {\rm 
 Can we actually compute the $Z$ and $V$ matrices of Theorem \ref{badsdp}, or the $U$ of Theorem \ref{goodsdp}?
 Regrettably,  
 we don't know how to do this in polynomial time either in the Turing model, or in the real number model of computing.
 However, we shall argue below that we can reduce this task to solving SDPs. 
 
 To start with the theoretical aspect of the reduction, we can find Z by running  a facial reduction
 algorithm 	\cite{BorWolk:81, WakiMura:12, Tuncel:11, Pataki:13}.  
 These algorithms must solve a sequence of SDPs in exact arithmetic.
 We can then verify whether ``Good condition 1" holds by solving the pair of SDPs 
 (\ref{redp})-(\ref{redd}). If it does hold, we can extract a $U$ matrix that satisfies it from an optimal solution of (\ref{redp}).  	
 If it does not, we can extract a $V$ certificate matrix that satisfies the ``Bad condition" from an optimal solution of the dual 
 (\ref{redd}).  
 
 In practice, heuristic and reasonably effective  implementations of facial reduction algorithms exist 
 \cite{permenter2014partial, zhu2019sieve}, and we may solve (\ref{redp})-(\ref{redd}) approximately, to deduce 
 that (\ref{p-sd}) is nearly badly or well behaved.  
 
 We mention here that the complexity of checking attainment and the existence of a positive gap in SDPs is unknown. 
 
 	} 
  \erem 

\section{When is the linear image of the semidefinite cone closed?}
\label{section-closed} 
We now address a question of independent interest in convex analysis/convex geometry:

\begin{center} 
	\framebox[4.5in]{Given a linear map, is the image of $\psd{n}$ under the map closed?}
\end{center}

This question fits in a much broader context. More generally, we can ask: when 
is the linear image of a closed convex set, say  $C,$  closed? Such closedness criteria are 
fundamental in convex analysis, and Chapter 9 in Rockafellar's classic text 
\cite{Rockafellar:70} is entirely dedicated to them. For more closedness  criteria see Chapter 2.3 in \cite{auslender2006asymptotic}, and 
for more recent work  on this subject
we refer to \cite{BausBor:99, BertTseng:07, BorweinMoors:09, BorweinMoors:10}. 
The latter paper shows that the 
set of linear maps under which the image of a closed convex cone is  {\em not} closed  is
small both in  measure and in category. 

 
 The closedness of the linear image of a closed convex cone
 ensures that a conic linear system is well-behaved 
 (in the same sense as (\ref{p-sd}));
 see e.g.,  \cite[Theorem 7.2]{Barvinok:2002}, 
 \cite[Theorem 2]{HenrionKorda:2014}, \cite[Lemma 2]{Pataki:17}. 
We studied criteria for the closedness of the 
linear image of a closed convex cone in \cite{Pataki:07}, 
and the results therein led to \cite{Pataki:17}, and to this paper. 
 
The special case $C = \psd{n}$ is interesting, since the semidefinite cone 
 is one of 
the simplest nonpolyhedral sets 
 whose geometry is well understood, see, e.g. 
\cite{BarCar:75, Pataki:00A} for a characterization of its faces. 
It turns out that the (non)closedness of the image of $\psd{n}$ admits 
simple combinatorial characterizations. 

We need some basic notation: for a set $S$ we define  its {\em frontier} $\front(S)$ as the difference between its closure and the set itself,
$$
\front(S) := \closure(S) \setminus S. 
$$
 
\begin{Example} \label{example-notclosed-2} 
	Define the map 
\begin{equation} \label{map1} 
	\sym{2} \ni Y   \rightarrow  (y_{11}, 2 y_{12}).  
	\end{equation} 
	The image of $\psd{2}$ -- shown on Figure \ref{figure-front1} in blue, and its frontier in red -- 
	is
		\begin{equation}
\{	(0,0) \} \cup \{ \, (\alpha, \beta): \, \alpha > 0 \, \}, 
	\end{equation}
	so it is not closed. For example,  $(0, 2)$ is in the frontier since 
		$
		(\epsilon, 2)$ is the image of the psd matrix 
		$$
		\bpx \eps & 1 \\ 1 & 1/\eps \epx \,\, 
		$$
		for all $\eps > 0, $ but no psd matrix is mapped to $(0, 2).$ 
			\begin{figure}[htp]
		\centering
		\includegraphics[width = 6cm]{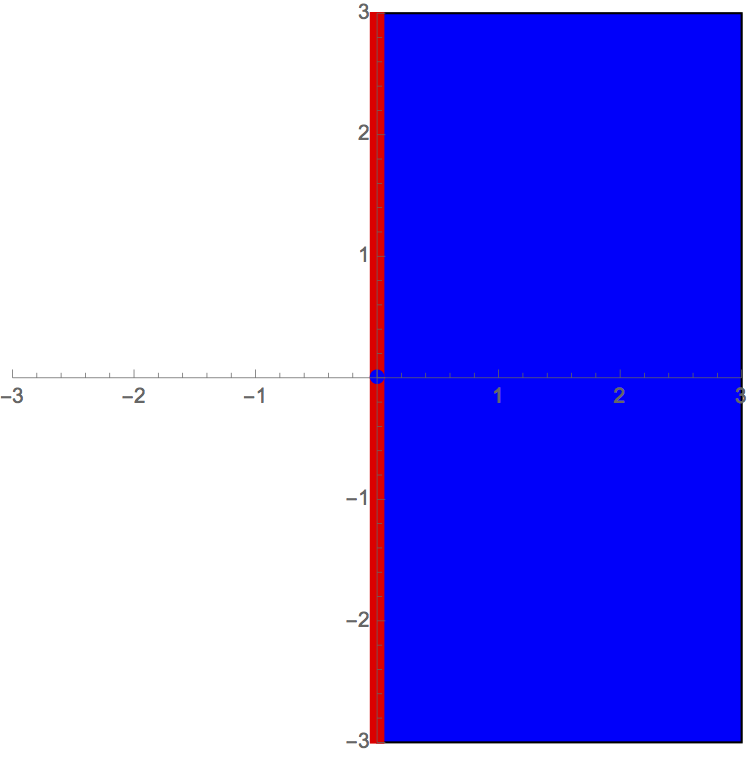}
			\caption{The image set is in blue, and its frontier is in red} 
			\label{figure-front1} 
	\end{figure}
\end{Example}

In more involved examples, however, the (non)closedness of the image 
 is much harder to check. 
\begin{Example} \label{example-notclosed-3} This example is based on 
	Example 6 in \cite{LiuPataki:17}. 
Define the linear map  
\begin{equation} \label{map2} 
\sym{3} \ni Y \rightarrow  ( 5 y_{11} + 4 y_{22} + 4 y_{13}, 3 y_{11} + 3 y_{22}  + 2 y_{13}, 2 y_{11} + 2 y_{22} + 2 y_{13}).
\end{equation} 
As we shall  see, the image of $\psd{3}$ is not closed, 
but verifying this by any ad hoc method seems very difficult.
\end{Example} 
For convenience, we shall  represent linear maps from 
$\sym{n}$ to $\rad{m}$ by matrices $A_1, \dots, A_m \in \sym{n}$ and write 
\begin{equation} \label{repr-A} 
\A(x) = \sum_{i=1}^m x_i A_i, \, {\rm and} \, \A^*(Y) \, = \, (A_1 \bullet Y, \dots, A_m \bullet Y).
\end{equation}
That is, we consider a linear map from $\sym{n}$ to $\rad{m}$ as the 
{\em adjoint} of a suitable linear map in the opposite direction, to better fit the framework of 
\cite{Pataki:07, Pataki:17}. 

The next proposition connects 
the closedness of the linear image of $\psd{n}$ and the
bad (or good) behavior of a homogeneous semidefinite system.
A  simple proof follows, e.g., 
	 from the classic separation theorem \cite[Theorem 1.1.1]{BorLewis:05}. 
	
\begin{Proposition} \label{proposition-closed} 
	Given a linear map $\A$ and its adjoint $\A^*$ as in (\ref{repr-A}), 
	the set $\A^*(\psd{n})$ is not closed if and only if the system 
	\begin{equation} \label{p-sd-h} \tag{\mbox{$P_{\mathit{SDH}}$}}
	\sum_{i=1}^m x_i A_i \preceq 0
	\end{equation}
	 is badly behaved. 
	 In particular, $c \in \front (\A^*(\psd{n}))$ if and only if the SDP 
	 \begin{equation}  
	 \begin{array}{rl} 
	 \sup &  \sum_{i=1}^m c_i x_i  \\
	 s.t.   & \sum_{i=1}^m x_i A_i \preceq 0 
	 \end{array}
	 \end{equation}
has optimal value zero, but its dual is infeasible. 
	 \qed 
\end{Proposition}

Thus, if (\ref{p-sd-h}) satisfies Assumption \ref{ass-slack}, then 
the characterizations of Theorems \ref{badsdp} and \ref{goodsdp} apply. 

More interestingly,  Corollary \ref{corollary-badlywell} and Proposition \ref{proposition-closed} together imply 
the following: 
\begin{Corollary} \label{corollary-closed} 
	Suppose $\A$ and $\A^*$ are represented as in 
	(\ref{repr-A}). Then $\A^*(\psd{n})$ is 
	\begin{enumerate}
		\item not closed if and only if the homogeneous 
		system (\ref{p-sd-h}) has a bad reformulation (of the form  (\ref{p-sd-bad}));
		\item closed if and only if the homogeneous 
		system (\ref{p-sd-h}) has a good reformulation (of the form (\ref{p-sd-good})).  
	\end{enumerate}
	\qed 
\end{Corollary}

We next illustrate Corollary \ref{corollary-closed} by continuing the previous examples. 
On the one hand, reformulating the map of  
Example \ref{example-notclosed-2} does not help either to verify nonclosedness of the image set, or to  
exhibit a vector in its frontier. Reformulating, however, does help a lot in Example \ref{example-notclosed-3}. 

\begin{Example} (Example \ref{example-notclosed-2}  continued) We can write the map 
	in (\ref{map1})  as 
$$
\sym{2} \ni Y \rightarrow \bpx  \bpx 1 & 0 \\ 0 & 0 \epx \bullet  Y, \, \bpx 0 & 1 \\ 1 & 0 \epx \bullet Y  \epx, 
$$
so the corresponding homogeneous semidefinite system is 
	$$
x_1 \bpx 1 & 0 \\ 0 & 0 \epx + x_2 \bpx 0 & 1 \\ 1 & 0 \epx \preceq 0,
$$
whose  bad reformulation is essentially the same: 
$$
x_1 \bpx 1 & 0 \\ 0 & 0 \epx + x_2 \bpx 0 & 1 \\ 1 & 0 \epx \preceq \bpx 1 & 0 \\ 0 & 0 \epx
$$
(we just replaced the the right hand side by the maximum rank slack). 
\end{Example} 

\begin{Example}
(Example \ref{example-notclosed-3}  continued) The homogeneous semidefinite system corresponding to the map in (\ref{map2}) is 
\begin{equation} \label{elaine-1} 
x_1 \bpx   5 &    0 &    2 \\
0 &    4 &    0 \\
2 &    0 &    0 \epx 
+ x_2 \bpx 3    & 0 &    1 \\
0 &    3 &    0 \\
1 &    0 &    0
\epx
+ x_3 \bpx 2   &  0 &    1 \\
0 &    2 &    0 \\
1 &    0  &    0
\epx \preceq 0. 
\end{equation}
Its bad reformulation  is 
\begin{equation} \label{jerry} 
x_1 \bpx 1 & 0 & 0 \\ 0 & 0 & 0 \\ 0 & 0 & 0 \epx 
+ x_2 \bpx  1 & 0 & 0 \\ 0 & 1 & 0 \\ 0 & 0 & 0     \epx 
+ x_3 \bpx 0 & 0 & 1 \\ 0 & 1 & 0 \\ 1 & 0 & 0 \epx  \preceq  \bpx  1 & 0 & 0 \\ 0 & 1 & 0 \\ 0 & 0 & 0     \epx.
\end{equation}

(\Newchange{How exactly did we obtain  (\ref{jerry})?}   To explain,  let us call the matrices $A_1, A_2,$  and 
$A_3$ on the left hand side in  (\ref{elaine-1}). 
Then (\ref{jerry}) is obtained by 
performing the
operations $A_2 = A_2 - A_3; A_1 = A_1 - 2 A_3; A_3 = A_3 - A_1 - A_2,$ then replacing the right hand  side by $A_2.$) 

Let $\A(x)$ be the left hand side in (\ref{elaine-1}) and 
$\A^{\prime}(x)$ the left hand side in (\ref{jerry}). 
Then  
$$
\A^{\prime *}(Y) = (y_{11}, y_{11} + y_{22}, y_{22} + 2 y_{13}),
$$
and a calculation shows 
(for details, see Example 6 in \cite{LiuPataki:17}) 
\begin{equation} \label{cla} 
\begin{array}{rcl} 
\closure (\A^{\prime *} \psd{3}) & = & \{ (\alpha, \beta, \gamma) \, : \, \beta \geq \alpha \geq 0 \, \},  \\ 
\front(\A^{\prime *} \psd{3}) & = & \{ (0, \beta, \gamma) \, | \, \beta \geq 0, \, \beta \neq \gamma  \, \}. 
\end{array} 
\end{equation}

The set  $\A^{\prime *}(\psd{3})$  is shown in Figure \ref{figure-jerry} in blue, and its frontier in red. 
Note that the blue diagonal segment on the red facet actually belongs 
to $\A^{\prime *}(\psd{3}).$
\begin{figure}[H]
	\centering
	\includegraphics[scale=0.5]{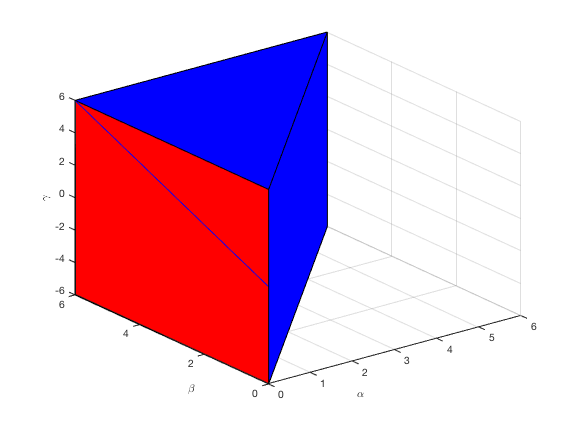} 
	\caption{The set $\A^{\prime *}(\psd{3})$ is in blue, and its frontier  in red} 
	\label{figure-jerry} 
\end{figure}
	The exact algebraic description of $\A^{\prime *}(\psd{3})$ (or of its closure and frontier) 	is 
	still not trivial  to find. 
However, its nonclosedness readily follows 
from Proposition 
\ref{proposition-closed} and Theorem 
\ref{badsdp}, since (\ref{jerry}) is 
badly behaved: we can choose $Z$ as the right hand side in
(\ref{jerry}) and $V$ as the coefficient matrix of $x_3.$ 

We can also quickly exhibit an element in   $\front(\A^{\prime *}(\psd{3})):$
the optimal value of the SDP 
\begin{equation} \nonumber 
\sup \, \{ \, x_3 \, | \, s.t. \, \A'(x) \preceq 0 \, \} 
\end{equation}
is  $0, \,$ but its dual is infeasible, 
hence   by Proposition \ref{proposition-closed} we deduce 
$$(0,0, 1) \in \front (\A^{\prime *} \psd{3}).$$ 
\end{Example}

	\begin{Remark} {\rm \label{remark-asymptote} 
\Newchange{We next connect our work to two other areas of convex analysis. The first area, 
asymptotes of convex sets, is classical; the second area, weak infeasibility in SDPs, 
 is more recent.} 

Let us define the distance of sets $S_1$ and $S_2$ 
as  
$$
\dist(S_1, S_2) \, := \, \inf \, \{ \, \norm{x_1 - x_2}  \, | \, x_1 \in S_1, \, x_2 \in S_2 \, \}.
$$
Let $ H := \, \{ \, Y \, | \, \A^*(Y) = c \, \}. \, $ Then by a standard argument  the following 
three statements are equivalent: 
\begin{enumerate}  
	\item \label{asymp-1}    $c \in \front(\A^*(\psd{n}));$
	\item  \label{asymp-2} $H \cap \psd{n} = \emptyset, \, {\rm and} \, \dist(H, \psd{n}) = 0;$
	\item \label{asymp-3} (\ref{sdp-d}) is infeasible, and its alternative system 
	\begin{equation} \label{eqn-alternative} 
	\begin{array}{rcl} 
	\sum_{i=1}^m c_i x_i & = & 1 \\
    \sum_{i=1}^m x_i A_i & \preceq  & 0 
    \end{array} 
    \end{equation} 
    is also infeasible. 
\end{enumerate} 

(The interested reader may want to work out the equivalences: for example,  one can use Theorem 11.4 in 
\cite{Rockafellar:70} which shows that two convex sets have a positive distance iff they can be separated in a strong sense.) 

 \Newchange{Note that whenever  (\ref{eqn-alternative}) happens to be {\em  feasible}, 
 it is an easy certificate that (\ref{sdp-d}) is {\em  infeasible}, as 
an argument analogous to proving weak duality shows 
that both cannot be feasible (hence the jargon  ``alternative system"). 
} 

Two terminologies are used to express the equivalent statements (\ref{asymp-1})-(\ref{asymp-3}) above. 

The first terminology says that $H$ is an {\em (affine) asymptote} of $\psd{n}.$ Asymptotes of convex sets were introduced  
in the classical paper \cite{klee1961asymptotes}. 
\Newchange{For example, 
	$$
	H = \{ \, Y \in \sym{2} \, | \, Y = \bpx 0 & 1 \\ 1 & y_{22} \epx \, {\rm \, for \, some} \, y_{22} \in \rad{}  \, \} 
	$$
is an asymptote of $\psd{2}:$ evidently $H$ and $\psd{2}$  do not intersect, but their 
distance is zero, since
$$
\bpx \eps & 1 \\ 1 & 1/\epsilon \epx \succeq 0 \,\,  {\rm for \, all} \,  \epsilon > 0. 
$$
} 
Alternatively, we can intersect  $\psd{2}$ with the hyperplane 
$\{ \, Y \in \sym{2} \, : y_{12}=y_{21} = 1 \, \}$ 
and check  that 
$\{ \, (0, y_{22}) \, : \, y_{22}  \in \rad{} \, \}$ is an asymptote of the resulting convex set (the area above a hyperbola). See 
the second part of Figure \ref{figure-hyperbola}. 

For more recent work  on asymptotes, see 
\cite{martinez2018minimization}, which shows that a convex set $C$ 
has an asymptote if and only if 
there is a quadratic function that is convex and lower bounded on $C$, but does not attain its infimum. 

The second terminology says that (\ref{sdp-d}) is {\em weakly infeasible}. 
Observe that when 
(\ref{sdp-p}) has finite optimal value and the dual (\ref{sdp-d}) is infeasible,  
it must be weakly infeasible. Indeed, suppose not; then the alternative system (\ref{eqn-alternative}) has a feasible solution $x$, and adding a large multiple of $x$ to a feasible solution of (\ref{sdp-p}) proves the latter is unbounded, which is a  contradiction.

	In more recent work 
	\cite{Lourenco:13} proved that a weakly infeasible SDP over $\psd{n}$  has a ``small" weakly infeasible subsystem of  dimension at most $n-1.$ 
	This result was generalized in Corollary 1 in \cite{LiuPataki:17} to conic linear programs, using 
	a  fundamental geometric 
	parameter of the underlying cone, 
	namely the length of the longest chain of faces.
	
}
\end{Remark}

\section{Discussion and conclusion} \label{section-conclusion}

We presented an elementary, in fact almost purely linear algebraic, 
proof of a combinatorial characterization of pathological semidefinite systems. 
En route, we showed how to transform semidefinite systems into normal 
forms to easily verify their pathological (or good) behavior. The normal forms also turned out to be useful for a related problem: they 
allow one to easily verify whether the linear image of $\psd{n}$ is closed. 

 We conclude with a discussion. 

\begin{itemize}
	
	\item 
	 As  we assumed throughout that (\ref{p-sd}) is feasible, 
	 we may ask:  does studying its bad behavior help us understand {\em all} pathologies in SDPs? 
	 
	 It certainly helps us understand many.  
	 In particular, it  
	 helps understand weak infeasibility, a pathology of infeasible SDPs: 
	 Remark \ref{remark-asymptote} and Proposition \ref{proposition-closed} show
	 that all $c$ that make (\ref{sdp-d}) weakly infeasible 
	 are suitable objective functions associated with badly behaved {\em homogeneous} (hence feasible) systems.
	 
	 However, we cannot yet distinguish among bad objective functions; 
	 for example, we cannot tell which $c \in \rad{m}$ gives a finite positive duality gap, 
	 and which gives the more benign pathology of zero duality gap coupled with unattained dual optimal value. 
	 
Since the interplay of semidefinite programming and algebraic geometry is a very active recent research area (some recent references are  \cite{Blekhetal:12, BhardRostalSanyal:15, nie2010algebraic, Vinzant:14}), it would be  interesting to 
connect
our results to  algebraic geometry.  

	\item Let us look again at the semidefinite systems in their normal forms  (\ref{p-sd-bad}) and (\ref{p-sd-good}) and note an interesting  feature they share. 
	They are both naturally split into two parts:  
\begin{itemize}
	\item a ``Slater part," namely the system 
	$\sum_{i=1}^k x_i F_i \preceq I_r, \,$ and 
	\item a  ``Redundant part," which corresponds to always zero variables 
	$x_{k+1}, \dots, x_m. \,$ 
	\end{itemize}
In 
(\ref{p-sd-bad}) the ``Redundant part" is responsible for the bad behavior. 

In 
(\ref{p-sd-good}) the ``Redundant part" is essentially linear: we can find the corresponding dual variable $Y_{22}$ by solving a system of equations, then doing a linesearch.

\item Here (and in \cite{Pataki:17}) we showed how  normal forms of semidefinite systems help  to 
verify their bad or good behavior. In more recent work, such normal forms turned out to be useful for other purposes: 
\begin{itemize}
	\item to verify the infeasibility of an SDP (see \cite{LiuPataki:15})  and 
	\item to verify the infeasibility and weak infeasibility of conic linear programs: see \cite{LiuPataki:17}. 
\end{itemize}

\item To construct  the normal forms,   
the bulk of the work is transforming the linear map
$$
\rad{m} \ni x \rightarrow \A(x) = \sum_{i=1}^m x_i A_i.
$$
Indeed, operations (\ref{exch})-(\ref{lambda}) of Definition \ref{definition-reform}  find an invertible linear  
map $M: \rad{m} \rightarrow \rad{m}$ so that $\A M$ is in an easier-to-handle form. 

Normal forms of linear maps are ubiquitous in linear algebra: see, for example, the row echelon form, or the eigenvector decomposition of 	a matrix. 	This work (as well as 
\cite{LiuPataki:15} and \cite{LiuPataki:17}) shows that 
they are also useful in a somewhat  unexpected area, the duality theory of conic linear 
programs. 

\end{itemize} 

{\bf Acknowledgement} I am grateful to the referees and the Area Editor for their detailed and helpful feedback; to Cedric Josz, Dan Molzahn,  and Hayato Waki for helpful discussions on SDP;  to Yuzixuan Zhu for her help with the figures; and to Yuzixuan Zhu and Alex Touzov for their careful reading of the paper. This research was supported by the National Science Foundation, award DMS-1817272.

\bibliography{mysdpMelody}

\end{document}